\documentclass[11pt]{article}
\usepackage{a4}
\usepackage[T1]{fontenc}
\usepackage[english]{babel}
\usepackage{latexsym}
\usepackage{amsmath}
\usepackage{amssymb}
\usepackage{amsthm}
\usepackage{mathrsfs}
\usepackage{color}
\usepackage{cite}
\usepackage{titling}

\newtheorem{teo}{Theorem}[section]
\newtheorem{defn}{Définition}[section]
\newtheorem{rem}{Remark}[section]

\newtheorem{lem}{Lemma}[section]

\usepackage{fancyhdr}

\fancypagestyle{plain}{\fancyhead{}}
\title{\textbf {Uniformly exponentially  stable approximations for a class of damped systems with unbounded feedbacks }}
\author{{Balegh Mohamed and Hajjej Zayd}\\ \textit{D\'epartement de Math\'ematiques, Facult\'e des Sciences de Monastir,}\\ \textit{Université de Monastir, 5019 Monastir, Tunisie.}\\ \textit{Unité de Recherche: Analyse et Contr\^ole des Equations aux Dérivées }\\ \textit{Partielles (Code : MESRS  05/UR/15-01)}.\\
\textit{E-mail address: hajjej.zayd@gmail.com}\\
\;\;\;\;\;\;\textit{E-mail address: Mohamedbalegh@gmail.com}}
\date{}
\begin{document}
\maketitle
 \abstract{}
 In this paper we study  time semi-discrete approximations of a class of exponentially stable infinite dimensional systems with unbounded feedbacks. It has recently been proved that for time semi-discrete systems, due to high frequency spurious components, the exponential decay property may be lost as the time step tends to zero. We prove that adding a suitable numerical viscosity term in the numerical scheme, one obtains approximations that are uniformly exponentially stable with respect to the discretization parameter.\\\\\




\section{Introduction}
Let $X$ and $Y$ be real Hilbert spaces ( $Y$ will be identified to its dual space) with norms denoted respectively by $\Vert .\Vert_X$ and $\Vert .\Vert_Y$.\\
Let $A: D(A)\to X$ be a skew-adjoint operator with  compact  resolvent and $B\in \mathfrak{L} (Y, D(A)')$, where $D(A)'$ is the dual space of $D(A)$ obtained by means of the inner product in $X$.\\
We consider the system described by
\begin{equation}\label{1}
\dot{z}(t)=Az-BB^*z,\;\;\;t\geq 0,\;\;\;z(0)=z_0\in X.
\end{equation}
Here and henceforth, a dot ($^.$ ) denotes differentiation with respect to time $t$. The element $z_0\in X$ is the initial state, and $z(t)$ is the state of the system.\\\\\\
Most of the linear equations modeling the damped vibrations of elastic structures can be written in the form (\ref{1}).
Some other relevant models, as the damped Schrodinger equations, fit in this setting as well.\\

We assume the following hypothesis introduced in \cite{4}:\\

(H)\\
If $\beta>0$ is fixed and $C_{\beta}=\{\lambda\in\mathbb{C} \;|\; Re \lambda=\beta\}$, the function
\begin{equation}\label{2}
\lambda\in\mathbb{C}_+=\{\lambda\in\mathbb{C}\; |\; Re \lambda>0\}\to H(\lambda)=B^*(\lambda I-A)^{-1}B\in \mathfrak{L}(Y)
\end{equation}
is bounded on $C_{\beta}$.\\

We define the energy of the solutions of system (\ref{1}) by:
\begin{equation}\label{3}
E(t)=\frac{1}{2}\Vert z(t)\Vert_X^2,\;\; t\geq 0,
\end{equation}
which satisfies
\begin{equation}\label{4}
\frac{d E}{dt}(t)=-\Vert B^*z(t)\Vert_Y^2,\;\;t\geq 0.
\end{equation}
In this paper, we assume that system (\ref{1}) is exponentially stable, that is there exist positive constants $\mu$ and $\nu$ such
that any solution of (\ref{1}) satisfies
\begin{equation}\label{5}
E(t)\leq \mu E(0)\exp(-\nu t),\;\; t\geq 0.
\end{equation}
Our goal is to develop a theory allowing to get, as a consequence of (\ref{5}), exponential stability results for time-discrete
systems.\\
We start considering the following natural time-discretization scheme for the continuous system (\ref{1}). For any $\Delta t$>0, we denote by $z^k$ the approximation of the solution $z$ of system (\ref{1}) at time $t_k=k\Delta t$, for $k\in\mathbb{N}$, and introduce
the following implicit midpoint time discretization of system (\ref{1}):
\begin{equation}\label{6}
\left \{
\begin{array}{lcr}
            \frac{z^{k+1}-z^k}{\Delta t} = A\left(\frac{z^k+z^{k+1}}{2}\right)-BB^*\left(\frac{z^k+z^{k+1}}{2}\right),\;k\in \mathbb{N}, & \\

             z^0=z_0.
\end{array}
\right.
\end{equation}
We define the discrete energy by:
\begin{equation}\label{7}
E^k=\frac{1}{2}\Vert z^k\Vert_X^2,\;\;\; k\in\mathbb{N},
\end{equation}
which satisfies the dissipation law
\begin{equation}
\frac{E^{k+1}-E^k}{\Delta t}=-\left\Vert B^*\left(\frac{z^k+z^{k+1}}{2}\right)\right\Vert_Y^2,\;\; k\in\mathbb{N}.
\end{equation}
It is well known that if the continuous system is exponentially  stable, the time-discrete ones do no more inherit of this property due to spurious high frequency modes ( see \cite{15}), that is we cannot expect in general to find positive constants $\mu_0$ and $\nu_0$ such that
\begin{equation}\label{1.1}
E^k\leq \mu_0 E^0\exp(-\nu_0k\Delta t),\;\;\;k\in\mathbb{N},
\end{equation}
holds for any solution of (\ref{6}) uniformly with respect to $\Delta t>0.$\\\\

Therefore, as in [12, 10, 13, 6], in order to get a uniform decay, it seems natural to add in system (\ref{6}) a suitable extra numerical viscosity term to damp these high-frequency spurious components. We obtain the new system:
\begin{equation}\label{8}
\left \{
\begin{array}{lcr}
            \frac{\tilde{z}^{k+1}-z^k}{\Delta t} = A\left(\frac{z^k+\tilde{z}^{k+1}}{2}\right)-BB^*\left(\frac{z^k+\tilde{z}^{k+1}}{2}\right),\;\;k\in \mathbb{N}, & \\\\

             \frac{z^{k+1}-\tilde{z}^{k+1}}{\Delta t}={(\Delta t)^2}{A^2}{z^{k+1}},\;\;\;k\in\mathbb{N},\\\\

             z^0=z_0.
\end{array}
\right.
\end{equation}
The energy of (\ref{8}), still defined by (\ref{7}), now satisfies:
\begin{equation}\label{9}
\left \{
\begin{array}{lcr}
            \tilde{E}^{k+1}=E^k - \Delta t \left\Vert B^*\left(\frac{z^k+\tilde{z}^{k+1}}{2}\right)\right\Vert_Y^2,\;\;\;k\in \mathbb{N}, \\\\

             E^{k+1}+(\Delta t)^3\Vert Az^{k+1}\Vert_X^2+\frac{(\Delta t)^6}{2}\Vert {A^2}{z}^{k+1}\Vert_X^2=\tilde{E}^{k+1},\;\;k\in \mathbb{N}.
 \end{array}
\right.
\end{equation}
Putting these identities together, we get:
\begin{equation*}
E^{k+1}+(\Delta t)^3\left\Vert Az^{k+1}\right\Vert_X^2+\frac{(\Delta t)^6}{2}\left\Vert {A^2}{z}^{k+1}\right\Vert_X^2+ \Delta t \left\Vert B^*\left(\frac{z^k+\tilde{z}^{k+1}}{2}\right)\right\Vert_Y^2 = E^k.
\end{equation*}
Summing this identities from $j=k_1$ to $j=k_2-1$, we obtain:\\

$\Vert z^{k_2}\Vert_X^2 + 2 \Delta t \displaystyle{\sum_{j=k_1}^{k_2-1}}\left\Vert B^*\left(\frac{z^j+\tilde{z}^{j+1}}{2}\right)\right\Vert_Y^2 + 2 \Delta t \displaystyle{\sum_{j=k_1}^{k_2-1}}(\Delta t)^2 \Vert Az^{j+1}\Vert_X^2 $\\

$ \hspace{5cm} + \Delta t \displaystyle{\sum_{j=k_1}^{k_2-1}}(\Delta t)^5 \Vert {A^2}{z^{j+1}}\Vert_X^2$

\begin{equation}\label{10}
=\Vert z^{k_1}\Vert_X^2,\;\;\; \forall k_1<k_2.
\end{equation}
The main result of this paper reads as follows:
\begin{teo}
Assume that system (\ref{1}) is exponentially stable, i.e. satisfies (\ref{5}) and the hypothesis (H) is verified.\\
Then there exist two positive constants $\mu_0$ and $\nu_0$ such that any solution of (\ref{8}) satisfies (\ref{1.1}) uniformly with respect to the discretization parameter $\Delta t>0$.
\end{teo}

Our strategy is based on the fact that the uniform exponential decay properties of the energy of systems (\ref{1}) and (\ref{8}) respectively are equivalent to uniform observability properties for the conservative system
\begin{equation}\label{11}
\dot{y}=Ay,\;\;\;t\in\mathbb{R},\;\;\;\;\;y(0)=y_0\in X,
\end{equation}
and its time semi-discrete viscous version:
\begin{equation}\label{12}
\left \{
\begin{array}{lcr}
            \frac{\tilde{u}^{k+1}-u^k}{\Delta t} = A\left(\frac{u^k+\tilde{u}^{k+1}}{2}\right),\;\;\;k\in \mathbb{N}, & \\\\

              \frac{u^{k+1}-\tilde{u}^{k+1}}{\Delta t}={(\Delta t)^2}{A^2}{u^{k+1}},\;\;\;k\in\mathbb{N},\\\\

            u^0=u_0.
 \end{array}
\right.
\end{equation}
At the continuous level the observability property consists in the existence of a time $T>0$ and a positive constant $k_T>0$ such that
\begin{equation}\label{13}
k_T\Vert y_0\Vert_X^2\leq \int_{0}^{T}\Vert B^* y(t)\Vert_Y^2\;dt,
\end{equation}
for every solution of (\ref{11}) (see \cite{4}).\\
A similar argument can be applied to the semi-discrete system (\ref{8}). Namely, the uniform exponential decay (\ref{1.1})
of the energy of solutions of (\ref{8}) is equivalent to the following observability inequality: there exist positive constants
$T$ and $c$ such that, for any $\Delta t>0$, every solution $u$ of (\ref{12}) satisfies:\\

$c\|u_0\|_X^2 \leq \Delta t \displaystyle { \sum_{k \Delta t \in \left[0,T\right]}}\|B^* u^k\|_Y^2+\Delta t \displaystyle { \sum_{k \Delta t \in \left[0,T\right]}}(\Delta t)^2\|Au^{k+1}\|_X^2$\\

\begin{equation}\label{14}
+\Delta t \displaystyle {\sum_{k \Delta t \in \left[0,T\right]}}(\Delta t)^5\|A^2 u^{k+1}\|_X^2.
\end{equation}
Our approach has common points with the result obtained in \cite{6} for feedbacks which are bounded in the energy
space. The main difference is that we replace the assumption of boundedness of $B$ by the assumption (H).\\
Let us mention the works  [13, 9], where boundary (that is $B$ is unbounded) stabilization  issues were discussed for space discrete of the  1-d wave equation.\\
In our knowledge, this paper is the first one providing exponential decay properties for time-discrete systems, when  the continuous setting has this property, in the case where $B$ is unbounded.\\
The outline of this paper is as follows.\\
In the second section, we give the background needed here. We recall some results on the observability of time-discrete conservative systems and prove (\ref{14}) in Section 3. Section 4 contains the proof of the main result. The last section is devoted to some applications.\\
In the following, to simplify the notation, $C$ and $c$ will denote a positive constants that may change from line to line, but don't depend on $\Delta t.$

\section{Some Background and Preliminaries}
In this section we give some background (without any proof) that we need in our present work ( for more details, see \cite{14}).\\
Throughout this section, $X$ is Hilbert space and $A:D(A)\to X$ be a densely defined operator with $\rho(A)\neq \emptyset$ ($\rho(A)$ is the resolvent set of $A$). We assume that $D(A)$ is endowed with the norm, $\Vert .\Vert_{D(A)} $, of the graph of $A$.\\
For every $\beta \in\rho(A)$, we define
$$ \Vert x\Vert_1=\Vert (\beta I-A)x\Vert_X,\;\;\;\forall x\in D(A).$$
The space $D(A)$ with this norm is a Hilbert space, denoted $X_1$. It is well known that $\Vert .\Vert_1$ is equivalent to $\Vert .\Vert_{D(A)}$.\\\\\\\\

We denote by $X_{-1}$ the completion of $X$ with respect to the norm
$$\Vert x\Vert_{-1}=\Vert (\beta I-A)^{-1}x\Vert_X \;\;\;\; \forall x\in X \;\text{and}\; \beta \in \rho(A).$$
Then $A\in\mathfrak{L}(X_1,X)$ and $A$ has a unique extension $\tilde{A}\in\mathfrak{L}(X,X_{-1})$.\\
Moreover,
\begin{equation}\label{15}
(\beta I-A)^{-1}\in\mathfrak{L}(X,X_1),\;\;\; (\beta I-\tilde{A})^{-1}\in\mathfrak{L}(X_{-1},X)
\end{equation}
and these two operators are unitary.\\

We recall also, if $A$ is maximal dissipative (for brevity m-dissipative)  then $(0,\infty)\subset \rho(A)$ and
\begin{equation}\label{16}
\Vert (\beta I-A)^{-1}\Vert_{\mathfrak{L}(X)}\leq \frac{1}{\beta}\;\;\;\forall \beta\in (0,\infty).
\end{equation}
When $A$ is skew-adjoint, we have  both $A$ and $-A$ are m-dissipative.\\
Now, we give the definition of a contraction.
\begin{defn}
A contraction is a bounded operator $C$, with the property that the norm $\Vert C\Vert \leq 1$. Note that the powers of a contraction have the same property, $\Vert C^n\Vert\leq 1$ for $n\in\mathbb{N}$.
\end{defn}
The following theorem (\cite{11}) gives a way of characterizing a generator of a semigroup of contraction.
\begin{teo}
Let $X$ be a Hilbert space, and let $A: D(A)\to X$ be a linear operator with dense domain. Then the following conditions are equivalent:\\

i) $A$ is the generator of a $C_0$ semigroup of contraction ,\\

ii) all $\lambda\in \mathbb{C}_+$ belong to $\rho(A)$, and $A_{\lambda}=(\bar{\lambda} I+A)(\lambda I-A)^{-1}$ is a contraction.

\end{teo}

Note that if $A$ is skew-adjoint, then it is the generator of a $C_0$ contraction semigroup.\\

Finally, we recall also that if $A$ is skew-adjoint, then we have $-A^2$ is a positive self adjoint operator and consequently $A^2$ is m-dissipative. \\

\section{Observability of time-discrete systems}

This section is organized as follows. First, we recall the results of \cite{5} on the observability of the time-discrete conservative system of (\ref{6}). Then, we give the proof of observability inequality (\ref{14}) which consists in the decomposition of the solution $u$ of (\ref{12}) into its low and high frequency parts, that we handle separately, as in \cite{6}.\\
\subsection{Some results on discrete observability}
We first need to introduce some notations.\\\\\\
Since $A$ is a skew-adjoint operator with compact resolvent, its spectrum is discrete and
$\sigma(A)=\{i\mu_j;\;j\in\mathbb{N}\}$, where $(\mu_j)_{j\in\mathbb{N}}$ is a sequence of real numbers such that $\vert \mu_j\vert \to\infty$ when $j\to\infty$. Set $(\phi_j)_{j\in\mathbb{N}}$ an orthonormal basis of eigenvectors of $A$ associated to the eigenvalues $(i\mu_j)_{j\in\mathbb{N}}$, that is
\begin{equation}\label{17}
A\phi_j=i\mu_j\phi_j.
\end{equation}
Moreover, define
\begin{equation}\label{18}
C_s(A) = span \{\phi_j:\;\text{the corresponding}\;i\mu_j\;\text{satisfies}\;\vert\mu_j\vert\leq s\}.
\end{equation}

The following theorem was proved in \cite{5}:
\begin{teo}
Assume that $B^*\in\mathfrak{L}(D(A),Y)$, that is
\begin{equation}\label{19}
\Vert B^*z\Vert_Y^2\leq C_B^2\Vert z\Vert_{D(A)}^2=C_B^2(\Vert Az\Vert_X^2+\Vert z\Vert_X^2),\;\;\forall z\in D(A),
\end{equation}
and that $A$ and $B^*$ satisfy the following hypothesis:
\begin{equation}\label{20}
\begin{cases}
\text{There exist constants}\;\; M, m>0\;\;\;\text{such that}\\\\
M^{2}\Vert(iwI-A)y\Vert_X^2+m^2\Vert B^*y\Vert_Y^2\geq \Vert y\Vert_X^2,\; \forall w \in\mathbb{R},\;y \in D(A).
\end{cases}
\end{equation}
\\
Then, for any $\delta>0$, there exists $T_{\delta}$ such that for any $T>T_{\delta}$, there exists a positive constant $k_{T,\delta}$, independent of $\Delta t$,  that depends only on $m, M, C_B, T$ and $\delta$, such that for $\Delta  t>0$ small enough, we have:
\begin{equation}\label{21}
 k_{T,\delta}\Vert y^{0}\Vert_X^2 \leq \Delta t \displaystyle { \sum_{k \Delta t \in \left[0,T \right]}
}\left\Vert B^*\left(\frac{y^k+y^{k+1}}{2}\right)\right\Vert_Y^2,\;\;  \forall y^0 \in C_{\delta/ \Delta t}(A),
\end{equation}
where $y^k$ is the solution of
\begin{equation}\label{22}
\frac{y^{k+1}-y^{k}}{\Delta t} = A\left(\frac{y^{k}+y^{k+1}}{2}\right),\;\;k\in \mathbb{N},\;\; y^0=y_0.
\end{equation}
\end{teo}

In the sequel, when there is no ambiguity, we will use the simplified notation $C_{\delta/\Delta t}$ instead of $C_{\delta/\Delta t}(A)$.\\
Hypothesis (\ref{20}) is the so-called Hautus test or resolvent estimate, which has been proved in \cite{8} to be equivalent to the continuous
observability inequality (\ref{13}) for the conservative system (\ref{11}) for suitable positive constants $T$ and $k_T$,  which
turns out to be equivalent to the exponential decay property (\ref{5}) for the continuous damped system (\ref{1}).\\

The following Lemma gives a resolvent estimate which was proved in \cite{6}:
\begin{lem}
Under the assumptions of Theorem 1.1, the resolvent estimate (\ref{20}) holds, with constants $m$ and $M$ that depend only on $\mu$ and $\nu$ given by (\ref{5}).
\end{lem}

Applying Theorem 3.1, for any $\delta>0$, choosing a time $T^*>T_{\delta}$ there exists a positive constant $k_{T^*,\delta}$ such that the inequality (\ref{21}) holds for any solution $y$ of (\ref{22}) with $y^0\in C_{\delta/\Delta t}$. In the sequel, we fix a positive number $\delta>0$ (for instance $\delta=1$), and $T^*=2T_{\delta}$.
\subsection{Uniform observability inequalities}
\begin{lem}
There exists a constant $c>0$ such that (\ref{14}) holds with $T=T^*$ for all solutions $u$ of (\ref{12}) uniformly with respect to $\Delta t$.
\end{lem}
{\bf Proof.}\;\; We decompose the solution $u$ of (\ref{12}) into its low and high frequency parts. To be more precise, we consider
\begin{equation}\label{24}
u_l=\pi_{\delta/\Delta t}u,\;\;\;\; u_h=(I-\pi_{\delta/\Delta t})u,
\end{equation}
where $\delta$ is the positive number that we have been chosen above, and $\pi_{\delta/\Delta t}$ is the orthogonal projection on $C_{\delta/\Delta t}$ defined in (\ref{18}). Here the notations $u_l$ and $u_h$ stand for the low and high frequency components, respectively.\\
Note that both $u_l$ and $u_h$ are solutions of (\ref{12}).\\
Besides, $u_h$ lies in the space $C_{\delta/ \Delta t}^\bot$, in which the following property holds:
\begin{equation}\label{25}
\Delta t \Vert Ay\Vert_X\geq \delta \Vert y\Vert_X,\;\;\;\;\forall y\in C_{\delta/ \Delta t}^\bot.
\end{equation}
{\bf The low frequencies.}\;\; We compare $u_l$ with $y_l$ solution of (\ref{22}) with initial data $y_l(0)=u_l(0)$. Set $w_l=u_l-y_l$. From (\ref{21}), we get:
$$k_{T^*,\delta}\Vert u_l^{0}\Vert_X^2\leq 2\Delta t \displaystyle { \sum_{k \Delta t \in \left[0,T^* \right]}
}\left\Vert B^*\left(\frac{u_l^k+\tilde{u}_l^{k+1}}{2}\right)\right\Vert_Y^2$$
\begin{equation}\label{26}
\hspace{4.1cm} + 2\Delta t \displaystyle { \sum_{k \Delta t \in \left[0,T^* \right]}
}\left\Vert B^*\left(\frac{w_l^k+\tilde{w}_l^{k+1}}{2}\right)\right\Vert_Y^2.
\end{equation}
The equation satisfied by $w_l$ is:
\begin{equation}\label{27}
\left \{
\begin{array}{lcr}
            \frac{\tilde{w}_l^{k+1}-w_l^k}{\Delta t} = A\left(\frac{w_l^k+\tilde{w}_l^{k+1}}{2}\right),\;\;k\in \mathbb{N},  \\\\

              \frac{w_l^{k+1}-\tilde{w}_l^{k+1}}{\Delta t}={(\Delta t)^2}{A^2}{u_l^{k+1}},\;\;k\in\mathbb{N},\\\\

            w_l^0=0.
 \end{array}
\right.
\end{equation}
It is easy to see that $\Vert \tilde{w}_l^{k+1}\Vert_X^2=\Vert w_l^k\Vert_X^2$. Besides, we have:\\

$\left\Vert B^*\left(\frac{w_l^k+\tilde{w}_l^{k+1}}{2}\right)\right\Vert_Y^2 \leq \frac{1}{2}\big\Vert B^*w_l^k\big\Vert_Y^2 + \frac{1}{2}\big\Vert B^*\tilde{w}_l^{k+1}\big\Vert_Y^2$.\\

Simple Calculations give:
$$B^*\tilde{w}_l^{k+1}=B^*(I-\frac{\Delta t}{2} A)^{-1}(I+\frac{\Delta t}{2} A)w_l^k.$$
Therefore, we have:\\
$$\big\Vert B^*\tilde{w}_l^{k+1}\big\Vert_Y^2  = \big\Vert B^*(I-\frac{\Delta t}{2} A)^{-1}(I+\frac{\Delta t}{2} A)w_l^k\big\Vert_Y^2$$

                                  $\hspace{3.5cm}  \leq  C_B^2 \big\Vert (I-\frac{\Delta t}{2} A)^{-1}(I+\frac{\Delta t}{2} A)w_l^k\big\Vert_{D(A)}^2$\\

                                  $\hspace{3.5cm} \leq  C \big\Vert (I-\frac{\Delta t}{2} A)^{-1}(I+\frac{\Delta t}{2} A)w_l^k\big\Vert_1^2.$\\\\\\

Finally, we get

                               $$ \big\Vert B^*\tilde{w}_l^{k+1}\big\Vert_Y^2\leq  C\big\Vert (I+\frac{\Delta t}{2} A)w_l^k\big\Vert_X^2$$

                                 $\hspace{4.5cm}\leq  C\big\Vert w_l^k\Vert_X^2 +c\frac{(\Delta t)^2}{4}\Vert Aw_l^k\big\Vert_X^2$\\

                                  $\hspace{4.5cm} \leq  C\big\Vert w_l^k\big\Vert_X^2,$

where we have used (\ref{15}) and the fact that $w_l^k\in C_{\delta/\Delta t}$.\\

In the same manner, by using (\ref{15}) and the fact that  $\tilde{w}_l^{k+1}\in C_{\delta/\Delta t}$ \\

and  $\Vert \tilde{w}_l^{k+1}\Vert_X^2=\Vert w_l^k\Vert_X^2$, we show that:\\

$$\big\Vert B^*w_l^k\big\Vert_Y^2 \leq c  \big\Vert w_l^k\big\Vert_X^2,$$

Consequently, we obtain:
\begin{equation}\label{28}
\left\Vert B^*\left(\frac{w_l^k+\tilde{w}_l^{k+1}}{2}\right)\right\Vert_Y^2 \leq c \left\Vert w_l^k\right\Vert_X^2.
\end{equation}
In addition, we have (see \cite{6})
\begin{equation}\label{29}
\big\Vert w_l^k\big\Vert_X^2 \leq \Delta t \sum_{j=1}^{k} (\Delta t)^2\big\Vert Au_l^j\big\Vert_X^2+\delta^2 \Delta t\sum_{j=0}^{k}\big\Vert w_l^j\big\Vert_X^2.
\end{equation}
Gr\"{o}onwal's Lemma applies and allows to deduce from (\ref{26}) and (\ref{28}), the existence of a positive $c$ independent of $\Delta t$ such that
\begin{equation*}
c\Vert u_l^{0}\Vert_X^2\leq \Delta t \displaystyle { \sum_{k \Delta t \in \left[0,T^* \right]}
}\left\Vert B^*\left(\frac{u_l^k+\tilde{u}_l^{k+1}}{2}\right)\right\Vert_Y^2 +\Delta t \displaystyle { \sum_{k \Delta t \in \left]0,T^*\right]}}(\Delta t)^2\|Au_l^k\|_X^2.
\end{equation*}

Besides,

$$\Delta t \displaystyle { \sum_{k \Delta t \in \left[0,T^* \right]}
}\left\Vert B^*\left(\frac{u_l^k+\tilde{u}_l^{k+1}}{2}\right)\right\Vert_Y^2\leq 2\Delta t \displaystyle { \sum_{k \Delta t \in \left[0,T^* \right]}
}\left\Vert B^*\left(\frac{u^k+\tilde{u}^{k+1}}{2}\right)\right\Vert_Y^2$$
$$ \hspace{6cm}+ 2\Delta t \displaystyle { \sum_{k \Delta t \in \left[0,T^* \right]}
}\left\Vert B^*\left(\frac{u_h^k+\tilde{u}_h^{k+1}}{2}\right)\right\Vert_Y^2.$$
Let us then estimate the last term in the right-hand side of the above inequality. We have

$$\left\Vert B^*\left(\frac{u_h^k+\tilde{u}_h^{k+1}}{2}\right)\right\Vert_Y^2 \leq \frac{1}{2}\big\Vert B^*u_h^k\big\Vert_Y^2 + \frac{1}{2}\big\Vert B^*\tilde{u}_h^{k+1}\big\Vert_Y^2.$$
But,
$$B^*u_h^k=B^*(I+\frac{\Delta t}{2} A)^{-1}(I-\frac{\Delta t}{2} A)\tilde{u}_h^{k+1}.$$\\\\

Then, we have

$$\big\Vert B^*u_h^k\big\Vert_Y^2  \leq  C \big\Vert (I+\frac{\Delta t}{2} A)^{-1}(I-\frac{\Delta t}{2} A)\tilde{u}_h^{k+1}\big\Vert_{D(A)}^2$$

                                $\hspace{2.8cm}   \leq  C \big\Vert (I-\frac{\Delta t}{2} A)\tilde{u}_h^{k+1}\big\Vert_X^2$\\

                                $\hspace{2.8cm}   \leq  2C \big\Vert \tilde{u}_h^{k+1}\big\Vert_X^2 +C \frac{(\Delta t)^2}{2}\big\Vert A\tilde{u}_h^{k+1}\big\Vert_X^2.$\\

Now, since $A\in\mathfrak{L}(X_1,X)$ and from (\ref{15}), we get

$$\big\Vert A\tilde{u}_h^{k+1}\big\Vert_X^2  =  \big\Vert A(I-\frac{\Delta t}{2} A)^{-1}(I+\frac{\Delta t}{2} A)u_h^k\big\Vert_X^2$$

                                           $\hspace{3.5cm}  \leq  C \big\Vert (I-\frac{\Delta t}{2} A)^{-1}(I+\frac{\Delta t}{2} A)u_h^k\big\Vert_{X_1}^2$\\

                                           $\hspace{3.5cm}   \leq  C  \big\Vert (I+\frac{\Delta t}{2} A)u_h^k\big\Vert_X^2$\\

                                           $\hspace{3.5cm}  \leq  2C \left(\big\Vert u_h^k\Vert_X^2 +\frac{(\Delta t)^2}{2} \big\Vert Au_h^k\big\Vert_X^2\right).$ \\

Consequently, for $\Delta t$ small enough, we obtain:

$$\big\Vert B^*u_h^k\big\Vert_Y^2\leq C \left(\big\Vert u_h^k\Vert_X^2 + (\Delta t)^2\big\Vert Au_h^k\big\Vert_X^2\right).$$

By the same way, we prove that:

$$\big\Vert B^*\tilde{u}_h^{k+1}\big\Vert_Y^2\leq C \left(\big\Vert u_h^k\Vert_X^2 + (\Delta t)^2\big\Vert Au_h^k\big\Vert_X^2\right),$$

and, sine $\tilde{u}_h^{k+1}$ and $u_h^k$ belong to  $C_{\delta/ \Delta t}^\bot$ for all $k$, we get from (\ref{25}) that

$\Delta t \displaystyle { \sum_{k \Delta t \in \left]0,T^* \right]}
}\left\Vert B^*\left(\frac{u_h^k+\tilde{u}_h^{k+1}}{2}\right)\right\Vert_Y^2 \leq C \Delta t \displaystyle {\sum_{k \Delta t \in \left]0,T^* \right]}
}\big\Vert u_h^k\big\Vert_X^2$ \\

$ \hspace{5.5cm} + C \Delta t \displaystyle { \sum_{k \Delta t \in \left]0,T^* \right]}
}(\Delta t)^2\big\Vert Au_h^k\big\Vert_X^2.$\\\\

$\hspace{5.5cm}\leq C \Delta t \displaystyle { \sum_{k \Delta t \in \left]0,T^* \right]}
}(\Delta t)^2\big\Vert Au_h^k\big\Vert_X^2.$\\

Therefore, we get:\\

$\Delta t \displaystyle { \sum_{k \Delta t \in \left[0,T^* \right]}
}\left\Vert B^*\left(\frac{u_h^k+\tilde{u}_h^{k+1}}{2}\right)\right\Vert_Y^2 \leq C \Delta t \displaystyle {\sum_{k \Delta t \in \left]0,T^* \right]}
}\big\Vert Au_h^k\big\Vert_X^2 $\\

$\hspace{5.5cm} + \Delta t\left\Vert B^*\left(\frac{u_h^0+\tilde{u}_h^{1}}{2}\right)\right\Vert_Y^2.$\\

Besides,\\

$\Delta t\left\Vert B^*\left(\frac{u_h^0+\tilde{u}_h^{1}}{2}\right)\right\Vert_Y^2 \leq 2 \Delta t\left\Vert B^*\left(\frac{u_l^0+\tilde{u}_l^{1}}{2}\right)\right\Vert_Y^2 + 2\Delta t\left\Vert B^*\left(\frac{u^0+\tilde{u}^{1}}{2}\right)\right\Vert_Y^2.$\\

Now, as (\ref{28}) we show that

$$\Delta t\left\Vert B^*\left(\frac{u_l^0+\tilde{u}_l^{1}}{2}\right)\right\Vert_Y^2\leq C\Delta t \big\Vert u_l^0\big\Vert_X^2.$$

Then,\\

$\Delta t \displaystyle { \sum_{k \Delta t \in \left[0,T^* \right]}
}\left\Vert B^*\left(\frac{u_h^k+\tilde{u}_h^{k+1}}{2}\right)\right\Vert_Y^2 \leq C \Delta t \displaystyle {\sum_{k \Delta t \in \left]0,T^* \right]}
}\big\Vert Au_h^k\big\Vert_X^2+C\Delta t \big\Vert u_l^0\big\Vert_X^2$ \\\\

$\hspace{6cm}+\Delta t\left\Vert B^*\left(\frac{u^0+\tilde{u}^{1}}{2}\right)\right\Vert_Y^2 .$\\

It follows that there exists $c>0$ independent of $\Delta t$ such that

$$c\Vert u_l^0\Vert_X^2\leq \Delta t \displaystyle { \sum_{k \Delta t \in \left[0,T^* \right]}
}\left\Vert B^*\left(\frac{u^k+\tilde{u}^{k+1}}{2}\right)\right\Vert_Y^2+\Delta t \displaystyle { \sum_{k \Delta t \in \left]0,T^*\right]}}(\Delta t)^2\|Au_l^k\|_X^2$$

\begin{equation}\label{30}
\hspace{3cm} +\Delta t \displaystyle { \sum_{k \Delta t \in \left]0,T^* \right]}
}(\Delta t)^2\big\Vert Au_h^k\big\Vert_X^2.
\end{equation}
{\bf The high frequencies.} We proceed as in \cite{6}, we get:\\

$C\Vert u_h^0\Vert_X^2\leq \Delta t \displaystyle { \sum_{k \Delta t \in \left]0,T^*\right]}}(\Delta t)^2\|Au_h^k\|_X^2$\\
\begin{equation}\label{31}
\hspace{5.5cm} +\Delta t \displaystyle { \sum_{k \Delta t \in \left]0,T^*\right]}}(\Delta t)^5\|A^2u_h^k\|_X^2.
\end{equation}
Combining (\ref{30}) and (\ref{31}) yields Lemma 3.2, since $u_h$ and $u_l$ lie in orthogonal spaces with respect to the scalar\\
 products $\langle .,.\rangle_X$ and $\langle A.,A.\rangle_X.$ \;\;$\square$
\section{Proof of Theorem 1.1.}
The proof of Theorem 1.1 will essentially rely on the following lemma:
\begin{lem}
Let $w$ be the solution of
\begin{equation*}
\left \{
\begin{array}{lcr}
            \frac{\tilde{w}^{k+1}-w^k}{\Delta t} = A\left(\frac{w^k+\tilde{w}^{k+1}}{2}\right)+Bv^k,\;\;\;k\in \mathbb{N}, & \\\\\

             \frac{w^{k+1}-\tilde{w}^{k+1}}{\Delta t}={(\Delta t)^2}{A^2}{w^{k+1}},\;\;\;\;\;k\in\mathbb{N},\\\\

             w^0=0,
\end{array}
\right.
\end{equation*}
where $v^k\in l^2(k\Delta t; Y)=\{v^k,\;\; \text{such that}\;\; \displaystyle{\sum_{k\Delta t\in[0,T^*]}}\Vert v^k\Vert_Y^2 <\infty\}$.\\

Let $T^*>0$ defined as before. There exists a positive constant $C$ independent of $\Delta t$  such that for all $0 < \Delta t < 1$, we have the following estimate

\begin{equation}\label{32}
\Delta t \displaystyle { \sum_{k \Delta t \in \left[0,T^* \right]}
}\left\Vert B^*\left(\frac{w^k+\tilde{w}^{k+1}}{2}\right)\right\Vert_Y^2 \leq C \Delta t \displaystyle { \sum_{k \Delta t \in \left[0,T^* \right]}
}\left\Vert v^k\right\Vert_Y^2,
\end{equation}
where $v^k \in l^2(k\Delta t; Y)$.
\end{lem}
{\bf Proof.}
We denote by:
\begin{equation}\label{321}
S_{k}= \frac{\widetilde{w}^{k+1}+w^{k}}{2}.
\end{equation}
From (\ref{30}) we get:
\begin{equation}\label{331}
S_{k}=(I-\frac{\Delta t}{2}A)^{-1}w^{k}+\frac{\Delta
t}{2}(I-\frac{\Delta t}{2}A)^{-1}B v^{k}.
\end{equation}
Besides,
\begin{equation}\label{341}
w^{k+1}=Lw^{k}+\Delta t RBv^{k},
\end{equation}
where:
\begin{equation}\label{351}
\begin{cases}
R=(I-(\Delta t)^{3}A^{2})^{-1}(I-\frac{\Delta t}{2}A)^{-1}\\\\
L=(I-(\Delta t)^{3}A^{2})^{-1}(I-\frac{\Delta
t}{2}A)^{-1}(I+\frac{\Delta t}{2}A).
\end{cases}
\end{equation}
Then, for $k\geq 1$:
 \begin{equation}\label{361}
 w^{k}=\Delta t\sum_{l=0}^{k-1} L^{k-1-l}RBv^{l}.
\end{equation}
Combining (\ref{331}), (\ref{361}) we deduce that:
\begin{equation}\label{371}
S_{k}=\frac{\Delta t}{2}(I-\frac{\Delta t}{2}A)^{-1}Bv^{k}+\Delta
t\sum_{l=0}^{k-1}(I-\frac{\Delta t}{2}A)^{-1}L^{k-1-l}RBv^{l}.
\end{equation}
Consequently,
\begin{equation*}
B^*S_{k}=\frac{\Delta t}{2}B^*(I-\frac{\Delta t}{2}A)^{-1}Bv^{k}+\Delta
t\sum_{l=0}^{k-1}B^*(I-\frac{\Delta t}{2}A)^{-1}L^{k-1-l}RBv^{l}.
\end{equation*}
Using Cauchy Schwarz inequality and the hypothesis (H), it follows that there exists $c >0$ independent of $\Delta t$ such
that\\\\
$\big\Vert B^{*}S_{k}\big\Vert_Y^2\leq c \big\Vert
v^{k}\big\Vert_Y^2+(\Delta t)^2 \displaystyle{\sum_{l=0}^{k-1}}
\big\Vert B^{*}(I-\frac{\Delta t}{2}A)^{-1}L^{k-1-l}RB
\big\Vert_{\mathfrak{L}(Y)}^2\displaystyle{\sum_{l=0}^{k-1}}\big\Vert
 v^{l} \big\Vert_Y^2.  $\\
Besides, we have:\\
$\big\Vert B^{*}(I-\frac{\Delta t}{2}A)^{-1}L^{k-1-l}RB
\big\Vert_{\mathfrak{L}(Y)}^2\leq \Vert B^*\Vert^2\Vert(I-\frac{\Delta t}{2}A)^{-1}\Vert^2\Vert L^{k-1-l}\Vert^2\Vert R\Vert^2\Vert B\Vert^2.$\\

Now, let us estimate  $\Vert L^{k-1-l}\Vert^2$. We have:\\

$\Vert L^{k-1-l}\Vert^2=\big\Vert [(I-(\Delta t)^{3}A^{2})^{-1}]^{k-1-l}[(I-\frac{\Delta
t}{2}A)^{-1}(I+\frac{\Delta t}{2}A)]^{k-1-l}\big\Vert^2$\\

$ \hspace{2cm} \leq\Vert (I+\frac{\Delta t}{2}A)(I-\frac{\Delta
t}{2}A)^{-1}\Vert^{2(k-1-l)}$\\

$\hspace{2cm}\leq 1,$\\

where we  used (\ref{16}) and Theorem 2.1.\\

Using again (\ref{16}) (since $A$ and $A^2$ are m-dissipative), we get\\

$\big\Vert B^{*}(I-\frac{\Delta t}{2}A)^{-1}L^{k-1-l}RB
\big\Vert_{\mathfrak{L}(Y)}^2\leq \Vert B\Vert^4.$\\

We deduce that\\

\begin{equation}\label{40}
\big\Vert B^{*}S_{k}\big\Vert_Y^2\leq c \big\Vert v^{k}
\big\Vert_Y^2 + (\Delta t)^{2}\big\Vert B \big\Vert^4
\frac{T}{\Delta t}\sum_{l=0}^{k-1}\big\Vert v^{l} \big\Vert_Y^2,
\end{equation}
which implies\\\\
\begin{equation}\label{41}
 \displaystyle{\sum_{k=1}^{\frac{T}{\Delta t}}} \big\Vert
B^{*}S_{k}\big\Vert_Y^2\leq c \displaystyle{\sum_{k=1}^{
\frac{T}{\Delta t}}} \big\Vert v^{k} \big\Vert_Y^2 +\big\Vert B
\big\Vert^4 T \Delta t \displaystyle{\sum_{k=1}^{\frac{T}{\Delta
t}}}\displaystyle{ \sum_{l=0}^{k-1}}\big\Vert v^{l} \big\Vert_Y^2.
\end{equation}
Since $k-1< \frac{T}{\Delta t}$, we obtain
\begin{equation}\label{42}
\displaystyle{\sum_{k=1}^{\frac{T}{\Delta t}}} \big\Vert
B^{*}S_{k}\big\Vert_Y^2\leq c \displaystyle{\sum_{k=1}^{
\frac{T}{\Delta t}}} \big\Vert v^{k} \big\Vert_Y^2 +\big\Vert B
\big\Vert^4 T \Delta t \displaystyle{\sum_{k=1}^{\frac{T}{\Delta
t}}}\displaystyle{ \sum_{l=0}^{\frac{T}{\Delta t}}}\big\Vert v^{l}
\big\Vert_Y^2
\end{equation}
\begin{equation}\label{43}
\hspace{3.3cm}\leq c \displaystyle{\sum_{k=1}^{ \frac{T}{\Delta t}}}
\big\Vert v^{k} \big\Vert_Y^2 +\big\Vert B \big\Vert^4 T \Delta t
\displaystyle{\sum_{l=0}^{\frac{T}{\Delta t}}}\displaystyle{
\sum_{k=1}^{\frac{T}{\Delta t}}}\big\Vert v^{l} \big\Vert_Y^2
\end{equation}
\begin{equation}\label{44}
\hspace{1cm}\leq c \displaystyle{\sum_{k=1}^{ \frac{T}{\Delta t}}}
\big\Vert v^{k} \big\Vert_Y^2 +\big\Vert B \big\Vert^4 T^{2}
\displaystyle{ \sum_{k=0}^{\frac{T}{\Delta t}}}\big\Vert v^{k}
\big\Vert_Y^2.
\end{equation}
Simple Calculations give:
\begin{equation}\label{45}
\big\Vert B^{*}S_{0}\big\Vert_Y^2\leq c \big\Vert v^{0}
\big\Vert_Y^2.
\end{equation}
Combining (\ref{44}), (\ref{45}) and the fact that $S_{k}=
\frac{\widetilde{w}^{k+1}+w^{k}}{2}$, we get the existence of
constant $C$ independent of $\Delta t$
such that\\
\begin{equation*}
\sum_{k=0}^{\frac{T}{\Delta t}}\left\Vert
B^{*}(\frac{w^{k}+\widetilde{w}^{k+1}}{2})\right\Vert_{Y}^{2}
\leqslant C\sum_{k=0}^{\frac{T}{\Delta t}}\big\Vert
v^{k}\big\Vert_Y^2.\;\;\;\;\square
\end{equation*}
Now, we give the proof of our main result.\\

{\bf Proof of  Theorem 1.1.}\;\;Here we follow the argument in \cite{7}.\\

We decompose $z$ solution of (\ref{8}) as $z=u+w$ where $u$ is the solution of the system (\ref{12}) with initial data $u^0=z^0$. Applying Lemma 3.2 to $u=z-w$, we get:\\

$ c \Vert z^0\Vert_X^2\leq 2\Bigg(\Delta t \displaystyle{\sum_{k\Delta t\in [0,T^*]}}\left\Vert B^*\left(\frac{z^k+\tilde{z}^{k+1}}{2}\right)\right\Vert_Y^2$\\
$ + \Delta t \displaystyle{\sum_{k\Delta t\in [0,T^*[}}(\Delta t)^2\left\Vert Az^{k+1}\right\Vert_X^2 +\Delta t\displaystyle{\sum_{k\Delta t\in [0,T^*[}}(\Delta t)^5\left\Vert {A^2}{z}^{k+1}\right\Vert_X^2 \Bigg )$\\
 $+ 2\Bigg( \Delta t \displaystyle{\sum_{k\Delta t\in [0,T^*]}}\left\Vert B^* \left(\frac{w^k+\tilde{w}^{k+1}}{2}\right)\right\Vert_Y^2+ \Delta t \displaystyle{\sum_{k\Delta t\in [0,T^*[}}(\Delta t)^2\left\Vert Aw^{k+1}\right\Vert_X^2$\\
\begin{equation}\label{33}
+\Delta t\displaystyle{\sum_{k\Delta t\in [0,T^*[}}(\Delta t)^5\left\Vert {A^2}{w}^{k+1}\right\Vert_X^2 \Bigg).
\end{equation}
Below, we bound the terms in the right-hand side of (\ref{33}) involving $w$ by the ones involving $z$.\\
The function $w$ satisfies:
\begin{equation}\label{34}
\left \{
\begin{array}{lcr}
            \frac{\tilde{w}^{k+1}-w^k}{\Delta t} = A\left(\frac{w^k+\tilde{w}^{k+1}}{2}\right)-BB^{\ast}\left(\frac{z^k+\tilde{z}^{k+1}}{2}\right),\;\;\;k\in \mathbb{N}, & \\\\\

             \frac{w^{k+1}-\tilde{w}^{k+1}}{\Delta t}={(\Delta t)^2}{A^2}{w^{k+1}},\;\;\;\;k\in\mathbb{N},\\\\

             w^0=0.
 \end{array}
\right.
\end{equation}
\\

By applying now Lemma 4.1 with $v^k =-B^*\left(\frac{z^k+\tilde{z}^{k+1}}{2}\right)$, we obtain that
\begin{equation}\label{35}
\displaystyle{\sum_{k\Delta t\in [0,T^*]}}\left\Vert B^*\left(\frac{w^k+\tilde{w}^{k+1}}{2}\right)\right\Vert_Y^2 \leq C \Delta t \displaystyle{\sum_{k\Delta t\in [0,T^*]}}\left\Vert B^*\left(\frac{z^k+\tilde{z}^{k+1}}{2}\right)\right\Vert_Y^2.
\end{equation}

Besides, we have (see \cite{6})
\begin{displaymath}
 \Vert w^k \Vert_X^2 + 2(\Delta t) \sum_{j=0}^{k-1}(\Delta t)^2 \Vert Aw^{j+1} \Vert_X^2 + (\Delta t)\sum_{j=0}^{k-1}(\Delta t)^5 \Vert {A^2}{w^{j+1}} \Vert_X^2
\end{displaymath}
\begin{equation}\label{36}
\leq \Delta t \sum_{j=0}^{k-1}\left(\left\Vert B^*\left(\frac{z^j+\tilde{z}^{j+1}}{2}\right)\right\Vert_Y^2+\left\Vert B^*\left(\frac{w^j+\tilde{w}^{j+1}}{2}\right)\right\Vert_Y^2\right).
\end{equation}
Combining (\ref{33}), (\ref{35}) and (\ref{36}) we get the existence of a constant $c$ independent of $\Delta t$ such that \\

$c \Vert z^0\Vert_X^2\leq 2\Delta t \displaystyle{\sum_{k\Delta t\in [0,T^*]}}\left\Vert B^*\left(\frac{z^k+\tilde{z}^{k+1}}{2}\right)\right\Vert_Y^2 + \Delta t \displaystyle{\sum_{k\Delta t\in [0,T^*[}}(\Delta t)^2\left\Vert Az^{k+1}\right\Vert_X^2$\\\\
 $$+\Delta t\displaystyle{\sum_{k\Delta t\in [0,T^*[}}(\Delta t)^5\left\Vert {A^2}{z}^{k+1}\right\Vert_X^2.$$
Finally, using the energy identity (\ref{10}), we get that
\[\Vert z^{T^{\ast}/ \Delta t}\Vert_X^2\leq (1-c)\Vert z^0\Vert_X^2.\]
The semi-group property then implies Theorem 1.1.\;\;\;\;\; $\square$

\begin{rem}
There are other possible discretizations schemes for system (\ref{1}) and other viscosity operators could have been chosen (see subsection 2.3 in \cite{6} for more details). Note that  the results given in the mentioned subsection still valid in our case but we replace the assumption of boundedness of $B$ by the assumption (H).
\end{rem}
\section{Applications}
\subsection{The wave equation}
We consider the following system
\begin{equation*}
\left\{
\begin{array}{lcr}

           u_{tt}(x,t)-u_{xx}(x,t)=0,\;\;x\in (0,\xi)\cup (\xi,1),\;t>0,\\\\

           u(0,t)=0,\;\;\;u_x(1,t)=0,\;\;t>0,\\\\

           u({\xi}_{-},t)=u({\xi}_{+},t),\;\;\;t>0,\\\\

           u_{x}( {\xi}_{-},t)-u_{x}({\xi}_{+},t) = -\alpha u_t(\xi,t),\;\;\;t>0,\\\\

            u(x,0)=u_0(x),\;\;u_t(x,0)=u_1(x),\;\;0<x<1,\end{array}
\right.
\end{equation*}

where $\xi\in (0,1)$ is a rational number with an irreducible fraction ($\xi=\frac{p}{q}$, where $p$ is odd), and $\alpha$ is a positive constant.\\

To show that this system  enters in the abstract setting of this paper, let us recall that it is equivalent to:\\

$\dot{Z}=AZ-BB^*Z,$\;\; with\; $Z=\left(\begin{array}{l}
u\\
v\\
\end{array}\right)$,\;\;$A=\left( \begin{array}{c}
0\;\;\;\;\;I \\
\partial_{xx}\;\;\;\;0
\end{array} \right).$\\
In this setting, $A$ is a skew-adjoint unbounded operator on the Hilbert space $X= V\times L^2(0,1)$, with domain $D(A)=H\times V$,\\
where
$$V=\{u\in H^1(0,1):\; u(0)=0\},$$
 and $$H=\{u\in H^2(0,1):\;u(0)=u_x(1)=0\}.$$
The operator $B$ is defined by:
$B:\mathbb{R}\to D(A)':\; k\to \left(\begin{array}{l}
0\\
\sqrt{\alpha}k\delta_{\xi}\\
\end{array}\right),$\;\; where $\delta_{\xi}$ is the Dirac mass concentrated in the point $\xi$.\\
It is well known that, in the case where $\xi=\frac{p}{q}$ ($p$ is odd), the energy of the system above decays exponentially, and the operators $A$ and $B$ defined above satisfy the assumption (H) (see \cite{2}).\\
Then, we introduce the following time semi-discrete approximation scheme:
\begin{equation}\label{38}
\left \{
\begin{array}{lcr}
            \frac{\tilde{Z}^{k+1}-Z^k}{\Delta t} = A\left(\frac{Z^k+\tilde{Z}^{k+1}}{2}\right)-BB^*\left(\frac{Z^k+\tilde{Z}^{k+1}}{2}\right),\;k\in \mathbb{N}, & \\\\

             \frac{Z^{k+1}-\tilde{Z}^{k+1}}{\Delta t}={(\Delta t)^2}{A^2}{Z^{k+1}} ,\;k\in\mathbb{N}, \\\\

             Z^0=(u_0,v_0).
\end{array}
\right.
\end{equation}
According to Theorem 1.1, we get:
\begin{teo}
There exist positive constants $\mu_0$ and $\nu_0$ such that any solution of (\ref{38}) satisfies (\ref{1.1}) uniformly with respect to the discretization parameter $\Delta t>0$.
\end{teo}
\subsection{One Euler-Bernoulli beam with interior damping}
We consider the following initial and boundary problem:
\begin{equation*}
\left\{
\begin{array}{lcr}
           u_{tt}(x,t)-u_{xxxx}(x,t)+\alpha u_t(\xi,t)\delta_{\xi}=0,\;\;0<x<1,\;t>0,\\\\

           u(0,t)=u_x(1,t)=u_{xx}(0,t)=u_{xxx}(1,t)=0,\;\;t>0,\\\\

           u(x,0)=u_0(x),\;\;u_t(x,0)=u_1(x),\end{array}
\right.
\end{equation*}
where $\xi\in (0,1)$ is defined as above, and $\alpha>0$. Hence it is written in the form (\ref{1}) with the following choices: Take $Z$ as above, and $X=H\times L^2(0,1)$. The operator $A$ defined by $$A=\left( \begin{array}{c}
0\;\;\;\;\;I \\
-\partial_{xxxx}\;\;\;\;0
\end{array} \right),$$
with domain $D(A)=V\times H$, where \\
$V=\{u\in H^4(0,1);\;u(0)=u_x(1)=u_{xx}(0)=u_{xxx}(1)=0\}$, and $H$ is defined as in the last subsection.\\
This operator is skew-adjoint on $X$. We now define the operator $B$ as:
$$B:\mathbb{R}\to D(A)':\;\;k\to \left(\begin{array}{l}
0\\
\sqrt{\alpha}k\delta_{\xi}\\
\end{array}\right).$$

The energy of the system above decays exponentially , and the hypothesis (H) was verified (see \cite{3}).\\
As an application of Theorem 1.1, we get:
\begin{teo}
The solutions of
\begin{equation*}
\left \{
\begin{array}{lcr}
            \frac{\tilde{Z}^{k+1}-Z^k}{\Delta t} = A\left(\frac{Z^k+\tilde{Z}^{k+1}}{2}\right)-BB^*\left(\frac{Z^k+\tilde{Z}^{k+1}}{2}\right),\;\;k\in \mathbb{N}, & \\\\

             \frac{Z^{k+1}-\tilde{Z}^{k+1}}{\Delta t}={(\Delta t)^2}{A^2}{Z^{k+1}},\;\;k\in\mathbb{N}, \\\\

             Z^0=(u_0,v_0),
\end{array}
\right.
\end{equation*}
are  exponentially uniformly decaying in the sense of (\ref{1.1}).
\end{teo}
\subsection{Dirichlet boundary stabilization of the wave equation}
Let $\Omega\subset\mathbb{R}^n,\; n\geq 2$ be an open bounded domain with a sufficiently smooth boundary  $\partial \Omega=\bar{\Gamma}_0\cup\bar{\Gamma}_1$, where $\Gamma_0$ and $\Gamma_1$ are disjoint parts of the boundary relatively open in $\partial \Omega$,\; $int(\Gamma_0)\neq \emptyset$. We consider the wave equation:
\begin{equation*}
\left \{
\begin{array}{lcr}
           u_{tt}-\Delta u=0,\;\;\Omega\times (0,+\infty),\\\\

           u=\frac{\partial}{\partial\nu}(Gu_t),\;\;\Gamma_0\times(0,+\infty),\\\\

           u=0,\;\;\;\Gamma_1\times (0,+\infty),\\\\

           u(x,0)=u_0(x),\;\;u_t(x,0)=u_1(x),\;\; \Omega,
\end{array}
\right.
\end{equation*}\\\\

where $\nu$ is the unit normal vector of $ \partial \Omega$ pointing towards the exterior of $\Omega$ and $G=(-\Delta)^{-1}: H^{-1}(\Omega)\to H_0^1(\Omega)$.\\
Denote: $A=\left( \begin{array}{c}
0\;\;\;\;\;I \\
\Delta\;\;\;\;0
\end{array} \right)$, with $D(A)=H_0^1(\Omega)\times L^2(\Omega)$.
In this setting, $A$ is a skew-adjoint unbounded operator on $X=L^2(\Omega)\times H^{-1}(\Omega)$. Moreover define
$$B\in\mathfrak{L}(L^2(\Gamma_0), D(A)'),$$
by $Bv=\left(\begin{array}{l}
0\\
-\Delta Dv
\end{array}\right)$,\;$\forall\; v\in L^2(\Gamma_0)$, where $D$ is the Dirichlet map i.e., $Df=g$ if and only if
\[\left \{
\begin{array}{lcr}
          \Delta g=0,\;\;\Omega,\\\\

          g=f,\;\Gamma_0,\;\;g=0,\;\Gamma_1.
\end{array}
\right.\]
The energy of the system above decays exponentially, and the hypothesis (H) was
verified (see \cite{1}).
Applying Theorem 1.1, we obtain:
\begin{teo}
The solutions of
\begin{equation*}
\left \{
\begin{array}{lcr}
            \frac{\tilde{Z}^{k+1}-Z^k}{\Delta t} = A\left(\frac{Z^k+\tilde{Z}^{k+1}}{2}\right)-BB^*\left(\frac{Z^k+\tilde{Z}^{k+1}}{2}\right),\;\;k\in \mathbb{N}, & \\\\

             \frac{Z^{k+1}-\tilde{Z}^{k+1}}{\Delta t}={(\Delta t)^2}{A^2}{Z^{k+1}},\;\;k\in\mathbb{N}, \\\\

             Z^0=(u_0,v_0),
\end{array}
\right.
\end{equation*}
are  exponentially uniformly decaying in the sense of (\ref{1.1}).
\end{teo}

\bibliographystyle{plain}

\end{document}